\documentclass[a4paper,12pt,twoside]{article}
\usepackage{xcolor}
\usepackage[latin1]{inputenc}
\usepackage[T1]{fontenc}
\usepackage{amssymb}
\usepackage{dochier}
\usepackage[final]{lpretex}
\usepackage{title}
\usepackage{a4}
\usepackage{amscd}
\usepackage{caption}
\input xy
\xyoption{all}

\EndOfDump

\def\DHpartformat#1{(#1) }
\def\DHrefpart#1{(\DHRefpart{#1})}
\LBsetstyleof{partno}{arabic}

\def\i {\item}

\let\define\def

  \def\C {{\mathbb C}}
 \def\E {{\mathbb E}} 
\def\GG {{\mathbb G}} \def\H{{\mathbb H}} \def\J {{\mathbb J}} 
  
  \def\P {{\mathbb P}} 
\def\Q {{\mathbb Q}}

\def\Z {{\mathbb Z}} 

\define \n {\mathbb N}
\define \q {\mathbb Q}
\define \PP {\mathbb P}

\def\sA {{\Cal A}}  \def\sC {{\Cal C}}
 \def\sE {{\Cal E}} \def\sF {{\Cal F}}
 \def\sH {{\Cal H}} 
  \def\sL {{\Cal L}}
\def\sM {{\Cal M}}  \def\sO {{\Cal O}}
 
 \def\sT {{\Cal T}} \def\sU {{\Cal U}}
  \def\sX {{\Cal X}}

\def\sZ {{\Cal Z}}
\define \cN {\Cal N}
\define \cf {\Cal F}
\define \cg {\Cal G}
\define \cE {\Cal E}
\define \ce {\Cal E}
\define \cc {\Cal C}
\define \cV {\Cal V}
\define \cA {\Cal A}
\define \cK {\Cal K}
\define \cO {\Cal O}
\define \cF {\Cal F}
\define \cn {\Cal N}
\define \cI {\Cal I}
\define \sP {\Cal P}
\define \sEll {\sE\ell\ell}
\define \sJE {\Cal{JE}}

\define \sGJE {\Cal{GJE}}
\define \sHyp {\mathcal{HYP}}
\def\a {\alpha} \def\b {\beta} \def\g {\gamma}

\def\s {\sigma}

\define \x {\xi}
\define \y {\eta}
\define \G {\Gamma}
\define \r {\rho}
\define \w {\omega}

\def \trho {\tilde {\rho}}

\define \tH {\widetilde H}
\define \tG {\widetilde{G}}
\define \tW {\widetilde W}
\define \tF {\widetilde F}
\define \tm {\tilde m}
\define \St {\widetilde S}
\define \Xt {\widetilde X}
\define \tS {\widetilde S}
\define \tpsi {\tilde \psi}
\define \tL {\widetilde L}
\define \tE {\widetilde E}
\define \tl {\tilde l}
\define \tA {\widetilde A}
\define \tom {\tilde\omega}
\define \tT {\widetilde T}
\define \tB {\widetilde B}
\define \tf {\tilde f}
\define \tsA {\widetilde{\sA}}

\define \tM {\widetilde M}
\define \tpsi {\widetilde{\psi}}
\define \trho {\widetilde{\rho}}
\define \tR {\widetilde R}
\define \tp {\widetilde p}
\define \tq {\widetilde q}
\define \tc {\tilde c}
\define \tsF {\widetilde {\sF}}
\define \tsM {\widetilde {\sM}}
\define \tii {\tilde i}
\define \tx {\tilde x}
\define \tg {\tilde g}
\define \tw {\tilde w}
\define \tz {\tilde z}
\define \ta {\widetilde\alpha}
\define \tb {\widetilde\beta}
\define \txi {\widetilde\xi}
\define \teta{\widetilde{\eta}}

\define \tsje{\widetilde{\sJE}}

\define \tLambda{\widetilde{\Lambda}}
\define \bD {\overline{D}}
\define \bG {\overline{G}}
\define \bI {\overline{I}}
\define \bK {\overline{K}}

\define\sje{\sJE}
\define\bsell{\bsEll}
\define \bV {\overline{V}}
\define \bX {\overline{X}}
\define \bY {\overline{Y}}
\define \bDelta {\overline{\Delta}}
\define \btau {\overline{\tau}}

\def\pd {\partial}

\def \Dx1 {\frac{\pd}{{\pd} x_1}}
\def \Dy1 {\frac{\pd}{{\pd} y_1}}
\def \Dz1 {\frac{\pd}{{\pd} z_1}}
\def \Dx2 {\frac{\pd}{{\pd} x_2}}
\def \Dy2 {\frac{\pd}{{\pd} y_2}}
\def \Dz2 {\frac{\pd}{{\pd} z_2}}

\def\q {\quad} 

\def\Mapdiagr#1{\nearrow\rlap{$\lower 5pt\vbox{{\hbox{$\mkern
-15mu\scriptstyle#1$}}}$}} 

\def\Mapdiagl#1{\llap{$\lower 5pt\vbox{{\hbox{$\scriptstyle#1\mkern
-15mu$}}}$}\searrow} 

\def\Mapswr#1{\swarrow\rlap{$\lower 5pt\vbox{{\hbox{$\mkern
-15mu\scriptstyle#1$}}}$}}              

\def\Mapnwl#1{\nwarrow\rlap{$\lower 5pt\vbox{{\hbox{$\mkern
-15mu\scriptstyle#1$}}}$}}

\def\i.e#1#2#3{\mathrel{\smash{\mathop{#2}\limits^{#1}_{#3}}}}

\def \inj {\hookrightarrow}

\define \Rhook {\hookrightarrow}

\def \half {\raise1pt\hbox{$\scriptstyle
        \frac{1}{2}\displaystyle$}}

\def \x{{\sl X}\llap{$\mkern -2mu {\scriptstyle -}$}}

\def \Bs {\operatorname{Bs}}

\def \Hom {\operatorname{Hom}}

\def \Proj {\operatorname{Proj}}

\def \Res {\operatorname{Res}}

\def \Bl {\operatorname{Bl}}
\def \Pic {\operatorname{Pic}}

\define \Kod {\operatorname{Kod}}
\define \dimension {\operatorname{dim}}
\define \codim {\operatorname{codim}}
\define \contr {\operatorname{contr}}
\define \rk {\operatorname{rank}}
\define \Im {\operatorname {Im}}
\define \Mor {\operatorname{Mor}}
\define \Cl {\operatorname{Cl}}
\define \Hilb {\operatorname{Hilb}}
\define \degree {\operatorname{deg}}
\define \mult {\operatorname{mult}}
\define \Aut {\operatorname{Aut}}
\define \NS {\operatorname{NS}}
\define \Gal {\operatorname{Gal}}
\define \ch {\operatorname{char}}
\define \Jac {\operatorname{Jac}}
\define \Km {\operatorname{Km}}
\define \Sec {\operatorname{Sec}}
\define \Stab {\operatorname{Stab}}
\define \Br {\operatorname{Br}}
\define \Inv {\operatorname {Inv}}
\define \tr {\operatorname{tr}}
\define \Frob {\operatorname{Frob}}
\define \Symn {\operatorname{Symm}^n}
\define \Ev {\sE^\vee}
\define \ordp {\operatorname{ord}_p}
\define \Supp {\operatorname{Supp}}
\define \Ann {\operatorname{Ann}}
\define \disc {\operatorname{disc}}
\define \lie {\operatorname{lie}}
\define \embdim {\operatorname{embdim}}

\def\Ram{\operatorname{Ram}}

\def\bsEll{\overline{\sE\ell\ell}}

\define\nbd{neighbourhood }

\define\Fil{\operatorname{Fil}}

\define\u0{\underline{0}}

\def\hod#1#2#3#4{\ensuremath{ if#30 H^{#2}({#1},{\cal O}_{#1}) \else 
 H^{#2}(#1,\Omega^{#3}\if\relax{#4}\relax_{#1}\else _{#1/#4}\fi)\fi}}

\begin{document}
\title[A differential equation]
{A non-linear differential equation for
  the {\MakeLowercase{de}} Rham cohomology
  of elliptic surfaces}
\author{N.I. Shepherd-Barron}
\address{King's College,
Strand,
London WC2R 2LS,
U.K.}
\email{Nicholas.Shepherd-Barron@kcl.ac.uk}
\maketitle
\begin{abstract}
  Suppose that $f:X\to C$ is a general
  Jacobian elliptic surface
  over $\C$. Write $\dim H^{1, 1}_{prim}(X)=N$;
  then $N$ is also the dimension of the stack $\sje$
  of Jacobian elliptic surfaces.
  There is, up to a sign,
  a natural orthonormal basis $(\eta_i)_{i\in [1, N]}$
  of $H^{1, 1}_{prim}(X)$ given by certain meromorphic
  $2$-forms $\eta_i$ of the second kind, one for each ramification
  point of the classifying morphism $\phi$ from $C$ to the
  stack $\bsell$ of generalized elliptic curves. A choice of local
  co-ordinate on $\bsell$ provides, via the branch locus of $\phi$,
  an {\'e}tale local co-ordinate system $(t_i)_{i\in [1, N]}$ on
  $\sje$.

  We prove here that
  some kind of truncation of the Gauss--Manin connection
  yields the system
  $$\{\partial_i H=(\partial_i \eta_i\wedge\eta_i)H\}_{i\in [1, N]}$$
  of non-linear pde satisfied by $H=[\eta_1,\ldots, \eta_N]$, where
  $\partial_i =\partial/\partial t_i$.
  Moreover, 
  $H$ can be interpreted as providing a multi-valued period map
  for these surfaces with values in the complex orthogonal group
  $O_N$, and we prove a generic
  infinitesimal Torelli theorem for this map.
  For rational elliptic surfaces this can be calculated explicitly.
\end{abstract}
AMS classification: 14C34, 32G20.
\begin{section}{Review of results on elliptic surfaces}
  Let $\sje=\sje_{h, q}$ denote the stack of complex
  Jacobian elliptic surfaces $f:X\to C$
  all of whose singular fibres are
  irreducible nodal curves; then $f$ is determined by
  the corresponding classifying morphism $\phi_f=\phi:C\to\bsell$
  where $\bsell\cong\P(4, 6)$ is the stack of generalized elliptic curves.
  This morphism $\phi$ is
  finite and is unramified over the locus
  $j=\infty$ (the cusp). Let $Z=\Ram_\phi=\sum (e_P-1)P$
  denote the ramification divisor, where $e_P$ is
  the ramification index of $\phi$ at $P\in C$.
  For any point $P\in C$, put $f^{-1}(P)=E_P$
  and let $M$ denote the line bundle of modular
  forms of weight $1$ on $\bsell$. One basic
  feature of $\sje$, some of whose details we recall below,
  is that we get an algebraic co-ordinate system
  on $\sje$ by taking
  the branch divisor $B=\phi_*Z$ in $\bsell$;
  $$\deg B=\deg Z=N=10 h+8(1-q)=\dim\sje_{h, q}.$$
\end{section}
\begin{section}{Moduli}
  In  \cite{SB1} we showed that
  the stack
  $\sje$ of Jacobian elliptic surfaces is smooth of the expected
  dimension $N=10h+8(1-q)$, provided that $10h>12(q-1)$,
  where $h$ is the geometric genus of $X$ and $q$ its irregularity,
  and the tangent bundle $T_{\sje}$
  is naturally a line
  bundle on the universal ramification divisor $\sZ$. That is,
  the tangent space $T_{\sje}(X)$ is given by
  $T_{\sje}(X)=\oplus_{P\in Z} L_P$ where $L_P$ is a skyscraper
  sheaf supported at $P$ and its length is $e_P-1$.

  We also know that $K_X\sim f^*(K_C+\phi^*M)$
  and $\deg\phi^*M=
  h-q+1=\chi(X,\sO_X)\ge 1$,
  so that $h\ge q$.
  In particular $\vert K_X\vert =f^*\vert K_C+\phi^*M\vert$
  and, by Riemann--Roch on $C$,
  has no base points if $h>q$. If $h=q$ then
  $\Bs\vert K_C+\phi^*M\vert$ is either empty or
  consists of a single point $Q$ and $\phi^*M\sim Q$.

  In fact
  the assumption that $10h>12(q-1)$ is unnecessary
  for the smoothness of $\sje$, although it does imply
  that the morphism $\sje\to\sM_q$ is smooth.

  \begin{theorem} $\sje$ is smooth of dimension $N$.
    \begin{proof} The exact sequence
      $$0\to T_C \to \phi^*T_{\bsell}\to \sL_Z\to 0,$$
      where $\sL_Z$ is a line bundle on $Z$, shows
      that $T_{\sje}(X)= H^0(Z,\sL_Z)$, which is of
      dimension $N$. On the other hand, for each point
      $P\in Z$ of ramification index $e_P=e$ we constructed
      in \cite{SB1} a variation of $X$ over
      an $(e-1)$-dimensional polydisc $\Delta^{e-1}$ whose tangent
      space is the part of $\sL$ that is supported at $P$. The
      product of these polydiscs is then the base of a smooth $N$-dimensional
      variation of $X$ whose tangent space is $H^0(Z,\sL_Z)$.
    \end{proof}
  \end{theorem}
  Moreover, by assigning to $f:X\to C$ the branch divisor of
  $\phi$ (that is, the image $B$ in $\bsell$
  of the ramification divisor $Z$ in $C$) we get a quasi-finite
  morphism $\sje\to\Hilb_N(\bsell)$. This defines an
  $\mathfrak S_N$-cover $\tsje$ of $\sje$ with a
  quasi-finite morphism $\tsje\to\bsell^N$, so that
  $\tsje$ has a natural system of local \emph{algebraic}
  co-ordinates and the polydisc $\Delta^{e-1}$ is, by
  its construction, defined
  by holding $N-(e-1)$ of these co-ordinates constant.
  In particular, therefore, $\Delta^{e-1}$ is contained
  in an algebraic substack of $\sje$.

  The forgetful morphism $\sje_{h, q}\to \sM_q$,
  the stack of curves of genus $q$, factors as
  $$\sje_{h, q}\buildrel\a\over\to \sU_{d, q}
  \buildrel\b\over\to\sM_q:
  (f:X\to C)\mapsto (C,\phi^*M)\mapsto (C).$$
  If $H^1(C,\phi^*T_{\bsell})=0$ (which is certainly the
  case if $2q-2<\deg\phi^*T_{\bsell}=10(1+h-q)$;
  that is, if $10 h> 8(q-1)$)
  then $\a$ is an open
  piece of a weighted projective space
  bundle whose fibre is 
  $\P(H^0(C, 4\phi^*M)\oplus H^0(C, 6\phi^*M))$
  (this parametrizes the coefficients of an
  equation for $X$ in Weierstrass form)
  and $\b:\sU_{d, q}\to\sM_q$ is a torsor
  under the universal Jacobian of degree $d$ that parametrizes
  the choice of line bundle $\phi^*M$ on $C$.
  (Here $d=\deg\phi^*M=1+h-q$.)

  Thus the description just given of the natural
  system of algebraic co-ordinates on $\sje_{h, q}$
  in terms of a branch divisor
  echoes Riemann's original conception (\cite{R}, p. 120)
  of moduli as natural algebraic co-ordinates (on the moduli space,
  although Riemann did not use that phrase)
  given by the points in $\P^1$ 
  of the branch divisor
  of a simply ramified morphism $\psi:C\to\P^1$.
  Rather, these co-ordinates must be normalized in some way,
  although Riemann did not indicate how this should be done.
  
  One might interpret
  this point of view to say that we first construct the moduli
  as the normalized branch points
  and then the moduli space is the object on which the moduli,
  which are multi-valued as far as distinguishing between
  curves is concerned,
  all become single-valued functions.
  If just one of the branch points is allowed to vary
  and the others held fixed then we get an algebraic
  curve in $\sM_q$ along which the derivative of the period map
  is of rank $1$; the difference
  from the picture in $\sje$ is that through each
  point of $\sM_q$ there is a positive-dimensional family
  of rank $1$ curves.

  For the rest of this paper we shall restrict attention to
  the open substack $\sje^{gen}$ of $\sje$ defined
  by the condition that $Z$ be reduced as well as
  disjoint from the cusp.
  Then \cite{SB1}, if
  $\vert K_X\vert$ has no base points\footnote{
    In \emph{loc. cit.} we stated this under the weaker assumption
    that $h>0$. The argument depended on the truth of the
    infinitesimal Torelli theorem, which holds if $\Bs\vert K_X\vert$ is
    empty \cite{S}, \cite{I}, \cite{Kl}.
    So the results of \cite{SB1} require the further assumption
    that $\Bs\vert K_X\vert$
    should be empty. However, if there
    are base points then things are more subtle; we shall
    discuss this below.}
  the derivative of the period map of the
  $1$-parameter variation $\sX\to\Delta$ to which
  reference has just been made defines,
  for every one of the points $P_1,\ldots, P_N\in\Ram_\phi$,
  a line $L_i$ in the vector space
  $\Hom(H^{2, 0}(X), H^{1, 1}_{prim}(X))$
  that is of rank $1$. More precisely,
  \begin{enumerate}
  \item the image of $L_i$ is a line $V_i$
    generated by the class $\eta_i$ of an explicit meromorphic
    $2$-form $\teta_i\in H^0(X,\Omega^2_X(2 E_P)) _{2k}$,
    taken modulo $H^{2, 0}(X)$,
    (the suffix $2k$ stands for ``of the second kind'') and
  \item the kernel of $L_P$ is $H^0(X, \Omega^2_X(-E_P))$.
    \item 
  The lines $L_i$ decompose $H^{1,1}_{prim}(X)$
  as the orthogonal
  direct sum of the lines $\C\eta_i$.
  The classes $\eta_i$
  form an orthogonal basis of $H^{1,1}_{prim}(X)$.
  This basis can then be rescaled to be orthonormal,
  up to a choice of sign.
  The proof of orthogonality in \cite{SB1} is of a variational
  nature and is only valid when $\vert K_X\vert$ has no base
  points; we shall show below that this orthogonality
  holds generally. We shall also
  give a description of the derivative of the period map
  even when there are base points.
\end{enumerate}
  Moreover, for any $1$-dimensional variation
  $\sX\to \Delta$ of $X$ the image 
  of the tangent space
  $T_{\Delta, 0}$ under the period map is a line in
  $\Hom(H^{2, 0}(X), H^{1, 1}_{prim}(X))$ that is of rank $1$
  if and only if $per_*(T_{\Delta, 0})$ equals one of the
  lines $L_i$.
  
  Then through any point of $\sje^{gen}$ there are $N$
  algebraic curves $\G_i,\ldots, \G_N$ that provide a co-ordinate frame at that
  point; each $\G_i$ is defined by fixing $N-1$ branch
  points on $\bsell$ and varying the remaining branch point.
  If the canonical linear system $\vert K_X\vert$ has no
  base points (see below for more about this assumption) then
  these curves are \emph{of rank $1$},
  in the sense that the derivative of the period map along any one of them
  has rank $1$ everywhere; this was proved in \cite{SB1}.

  \begin{proposition} Any holomorphic arc in $\sje^{gen}$ that is
  of rank $1$ is contained in one
  of these algebraic curves.
  \begin{proof} Any choice of algebraic local co-ordinate on $\bsell$
    provides algebraic co-ordinates $t_1,\ldots, t_N$ on $\sje^{gen}$ at a point
    $X$ of $\sje^{gen}$.
    There is then an orthogonal basis $(\eta_i)$ of
    $H^{1, 1}_{prim}(X)$ given by $\eta_i=\partial_{t_i}\omega$
    where $\omega\in H^{2, 0}(X)$ does not vanish along the
    fibres $E_P$ for $P\in Z$. (We can then rescale the $\eta_i$
    to make them orthonormal.) Under the period map
    the image of $\partial_{t_i}$ is a linear map whose
    kernel is $H^0(X,\Omega^2_X(-E_P))$ and whose
    image is in the line $\C\eta_i$.

    Suppose that $\Delta$ is a holomorphic
    arc in $\sje^{gen}$ through the point that is $X$
    that is of rank $1$. Say $s$ is a co-ordinate on $\Delta$.
    Then the image of $\partial_s$ under the
    period map is a homomorphism
    $H^{2, 0}(X)\to H^{1, 1}_{prim}(X)$ which is a linear
    combination of the images of the vectors $\partial_{t_i}$.
    Since it is of rank
    $1$ it is, from the description just given of the image
    of $\partial_{t_i}$, proportional to some $\partial_{t_i}$.
    Then, inside the tangent
    bundle to $\sje^{gen}$, $\partial_s$ is everywhere
    proportional to $\partial_{t_i}$, which means that
    $\Delta$ is contained in $\G_i$.
  \end{proof}
\end{proposition}
\end{section}

\begin{section}{A period matrix for the $(1, 1)$ classes}\label{3}
  As already remarked,
  any choice of local co-ordinate on $\bsell$
  gives an unordered {\'e}tale local co-ordinate system $\{t_1,\ldots, t_N\}$
  on $\sje^{gen}$ via the {\'e}tale
  morphism $\sje^{gen}\to\Hilb_N(\bsell)$,
  the stack that classifies the length $N$ closed subschemes
  of $\bsell$,
  which maps the surface $X\to C$
  to the branch locus in $\bsell$ of the classifying
  morphism $\phi: C\to \bsell$.
  
  Define $\sje^{gen, *}\to\sje^{gen}$ to be the
  ${\mathfrak S}_N$-cover induced by $\bsell^N\to\Hilb_N(\bsell)$.
  Then $(t_1, \ldots, t_N)$ is an ordered local co-ordinate
  system on $\sje^{gen, *}$. (We shall see below that
  if $h$ is not too small compared to $q$ then
  $\sje$ and $\sje^{gen, *}$ are irreducible.)

  \begin{lemma}
    
    \part[i] There is a unique finite $\Z/2$-cover
    $\bsell_+\to\bsell$ that is ramified over the cusp.

    \part[ii] $\bsell_+$ is isomorphic to the closed substack
    of $\P(4, 6, 6)$ defined by the equation
    $$g_4^3 -27 g_6^2=h_6^2.$$

    \begin{proof} \DHrefpart{i}
      The existence and uniqueness of $\bsell_+$
      is an immediate
      consequence of the description of $\bsell$
      as the modular curve $X(1)$, the unique normal compactification of
      the stack $\mathfrak H/SL_2(\Z)$.
      The cover $\bsell_+$ is associated to the
      unique index $2$ subgroup of $SL_2(\Z)$;
      this subgroup is the normal closure of the group
      generated by the square of the standard upper
      triangular unipotent matrix.

      \DHrefpart{ii} The
      description of $\bsell_+$ in terms of equations
      follows from the formula for the discriminant
      of an elliptic curve in Weierstrass form.
    \end{proof}
  \end{lemma}  

  Let $G_N=(\Z/2)^N\rtimes \mathfrak S_N$
  denote the group of signed permutations of $N$ letters.
  Then $\bsell_+^N\to \Hilb_N(\bsell)$ is a $G_N$-cover.
  Let $\tsje^{gen}\to\sje^{gen}$ denote the induced
  $G_N$-cover.
  \begin{lemma}\label{monodr}
    The rank $N$ vector bundle $\sH^{1, 1}_{prim}$
  on $\sje$ has a non-degenerate inner product
  and the pull back of $\sH^{1, 1}_{prim}$
  to $\tsje^{gen}$ is equipped with an orthonormal basis
  $(\eta_i)$ indexed by the ramification divisor.
  \begin{proof} There certainly is a $G_N$-cover
    $\sje^{\dagger}\to\sje^{gen}$ over which the pull back of
    $\sH^{1, 1}_{prim}$ has the structure described in the lemma.

    To see that $\sje^{\dagger}=\tsje^{gen}$,
    consider a $1$-parameter
    family $\sX\to\Delta$ in $\sje$ which is generic
    subject to the condition that 
    the ramification divisor $Z$ should specialize so as
    to remain reduced but meet the cuspidal locus.
    The closed fibre $X_0$ acquires a node
    while the total space $\sX$ is smooth.
    Then the monodromy on $H^{1, 1}_{prim}(X_t),$
    where $t\ne 0$,
    is a reflexion $\s$. It acts trivially on $\eta_j$
    if the index $j$ is attached to a point of $Z$
    that is not cuspidal and so
    $\s(\eta_i)=-\eta_i$ otherwise.
    It follows that $\sje^{\dagger}=\tsje^{gen}$.
  \end{proof}
\end{lemma}
  The line
  $\C\eta_i$ is the image of $H^0(X,\Omega^2_X)$
  under $\pd_i=\pd/\pd t_i$, provided that $\vert K_X\vert$
  has no base points (which is always the case if
  $h >q$) and the kernel of $\pd_i$ is
  $H^0(X, \Omega^2_X(-E_i))$. Write
  $$H=[\eta_1,\ldots, \eta_N].$$
  
  We next show that,
  after making a certain choice of cohomology
  classes on $X$, we can regard each $\eta_i$
  as a column vector of complex numbers
  and $H$ as an orthogonal complex $N\times N$ matrix that
  is a holomorphic function of the co-ordinates
  $t_1,\ldots, t_N$. 

  The primitive cohomology $H^2(X,\Z)_{prim}$ of a surface
  $X$ in $\sje^{gen}$ is even, unimodular and of signature
  $(2h, N)$, and so is isometric to the unique even unimodular
  lattice $\tLambda:=U^{2h}\perp E_8(-1)^{1+h-q}$
  of signature $(2h, 10h+8(1-q))$,
  where $U$ is the hyperbolic plane and $E_8$ is the
  root lattice of type $E_8$.

  Say that a sublattice $L$ of $\tLambda$ is \emph{good for $X$}
  if it is isotropic of rank $h$ and there is no line $\ell=\C\omega
  \subseteq H^{2, 0}(X)$ such that $L\subseteq \ell^\perp$.

  \begin{lemma}\label{surj} There exists a sublattice $L$ of $\tLambda$
    that is good for $X$.
    \begin{proof} Suppose not. The $O_\Lambda(\Z)$-orbit
      of $L$ is Zariski-dense in the Grassmannian
      $Gr$ of isotropic $h$-planes in $\tLambda_\C$, and so
      for every $V\in Gr$ there is a line
      $\ell\subset H^{2, 0}(X)$ such that $V\subseteq \ell^\perp$.
      An elementary argument using the incidence
      variety of such pairs $(V, \ell)$ gives a contradiction.
    \end{proof}
  \end{lemma}

  Pick an isotropic sublattice $L$ of $\tLambda$ whose
  rank is $h$ and put $\Lambda =L^\perp/L$, so that
  $\Lambda \cong U^h\perp E_8(-1)^{1+h-q}$,
  the unique even unimodular lattice of signature
  $(h, 9h+8(1-q))$. (Recall
  that any two such sublattices $L$ are equivalent
  under the orthogonal group $O_{\tLambda}(\Z)$.) Let $\J\E\to\tsje$
  be the holomorphic stack
  of pairs $(X, \psi)$ where $X$ is a surface in $\tsje$
  and $\psi:\tLambda\to H^2(X,\Z)_{prim}$ is an isometry
  and let $\J\E_L$ be the open substack of $\J\E$ defined by the condition
  that $\psi(L)$ is good for $X$.

  Observe that, since the complement of $\J\E_L$ is analytic in $\J\E$,
  for any pair $(X, \psi)$ in $\J\E^{gen}$
  a versal deformation $(\sX,\psi)\to B$ of $(X,\psi)$ has the property
  that $(\sX(t),\psi)$ lies in $\J\E_L^{gen}$.

  If $(X,\psi)\in\J\E_L^{gen}$
  then there is, for each $i$, a unique meromomorphic
  $2$-form $\teta_i\in H^0(X, \Omega^2(2E_i))_{2k}$
  such that $\eta_i=\teta_i$ modulo $H^{2, 0}(X)$
  and $\int_a\teta_i=0$ for all $a\in L$.

  Choose an arbitrary splitting of the surjection
  $L^\perp\to\Lambda$ and lift the elements of a chosen basis
  of $\tLambda/L^\perp$ to $b_1,\ldots, b_h$ of $\tLambda$.

  We can write $\teta_j=\sum_k\a_{j k}a_k+\lambda_j$
  with $\lambda_j\in \Lambda_\C$ and then
  $$\delta_{i j}=(\eta_i,\eta_j)=(\lambda_i,\lambda_j)=(\teta_i,\teta_j).$$
  Choose cycles $c_1,\ldots, c_N$ that
  are dual to the $\Z$-basis of $\Lambda$ that is specified by $\psi$;
  then $\lambda_j$ is the column vector
  $${}^t\bigg[\int_{c_1}\teta_j,\ldots, \int_{c_N}\teta_j\bigg].$$
  Each $\int_{c_i}\teta_j$ depends only on $i$ and $\eta_j$,
  so that we can write $\int_{c_i}\teta_j=\int_{c_i}\eta_j$.
  That is, we can regard each $\eta_j$ as a column vector,
  that is, as an element of
  $\C^N$, and $H$ as an $N\times N$ matrix
  $H=[\eta_1,\ldots, \eta_N]=[\lambda_1,\ldots,\lambda_N]$.

  Let $\sT_\Lambda$ denote the torsor under the orthogonal group
  $O_N$ that consists of
  $N\times N$ matrices whose
  columns are orthonormal with respect to the pairing on
  $\Lambda$. Then
  there is a holomorphic matrix-valued morphism
  $$H_L:\J\E_L^{gen}\to \sT_\Lambda$$
  defined by $H$.
  
  Recall that if $x, y$
  are elements of an inner product space then $x\wedge y$
  acts on a third element $z$ by
  $$(x\wedge y)z=(y, z)x - (x, z) y,$$
  the contraction of $x\wedge y$ against $z$.
  We extend this to an action of $x\wedge y$
  on a matrix $H=[z_1,\ldots, z_N]$ column by column.
  
  We can truncate the implicit linear system of pde on $H^2_{prim}$
  that is given by the Gauss--Manin connection
  to get an explicit non-linear system of pde that constrains
  $H^{1, 1}_{prim}$. In what follows we regard $H_L$ as a
  matrix-valued function of the local co-ordinates
  $t_1,\ldots, t_N$ on $\J\E_L^{gen}$ induced from
  the algebraic local co-ordinates $t_1,\ldots, t_N$ on
  $\tsje^{gen}$ and $\partial_i=\partial/\partial t_i$.
  
  \begin{theorem}\label{pde}
    $H_L:\J\E_L^{gen}\to \sT_\Lambda$ satisfies the system
    $$\{\pd_i H_L=(\pd_i\eta_i\wedge\eta_i)H_L\}_{i\in [1, N]}$$
    of non-linear pde.
    \begin{proof} We prove this column by column.
      
      We showed in \cite{SB1} that
      $\pd_i\eta_j$ is proportional to $\eta_i$ modulo $H^{2, 0}$
      when $i\ne j$. Since the $\eta_i$ are normalized
      by the constraint $\int_{a_k}\eta_l=0$ for all $k, l$
      there is then is a holomorphic
      function $f_{i j}$ on $\J\E_L^{gen}$, if $i\ne j$, such that
      $$\pd_i\eta_j=f_{i j}\eta_i.$$ Taking inner products
      shows that
      $(\pd_i\eta_j, \eta_i)=f_{i j},$
      so that
      $$\pd_i\eta_j=(\pd_i\eta_j, \eta_i)\eta_i.$$
      On the other hand
      $$(\pd_i\eta_i\wedge\eta_i)\eta_j=
      (\eta_i, \eta_j)\pd_i\eta_i-(\pd_i\eta_i, \eta_j)\eta_i
      =-(\pd_i\eta_i, \eta_j)\eta_i.$$
      But $(\eta_i, \eta_j)=0$, which differentiates to
      $$(\pd_i\eta_i, \eta_j)+ (\eta_i,\pd_i\eta_j)=0,$$
      and then
      $$(\pd_i\eta_i\wedge\eta_i)\eta_j=
      (\eta_i,\pd_i\eta_j)\eta_i=\pd_i\eta_j.$$
      So the equation holds when $i\ne j$.

      When $i=j$ the situation might appear to be more complex,
      because $\pd_i\eta_i$ merely exists as an element
      of $\C^N$: it no longer lies in $\sF il^1$.
      However, in this case
      the equation merely states that
      $\pd_i\eta_i =\pd_i\eta_i,$
      so that it holds in full.
    \end{proof}
  \end{theorem}
  We shall refer to this as the ecliptic equation
  since $\pd_i\eta_i\wedge\eta_i$ is the ecliptic of
  the particle $\eta_i$,
  the $2$-plane
  defined by its position
  and its velocity with respect to $t_i$, or the
  $2$-plane in which it is instantaneously moving.

  \begin{theorem} $H_L$ extends to a holomorphic
    map $\J\E^{gen}\to\sT$ that satisfies the ecliptic equation.
    \begin{proof} Given $(X, \psi)\in\J\E^{gen}$
      consider a versal deformation $\sX\to \sC\to B$
      and let $\sE_i\to B$ denote the $i$th curve in $\sX$
      where the classifying map $\sC\to\bsell\times B$
      is ramified.
      
      Say that the normalized form
      $\teta_i(t)\in H^0(\sX(t),\Omega^2(2\sE_i(t)))_{2k}$
      specializes to $\teta_i(0)$; then
      \begin{enumerate}
      \item $\teta_i(0)\in H^0(\sX, \Omega^2(2\sE_i(0)))_{2k}$,
      \item      $\int_a\teta_i(0)=0$ for all $a\in L$ and
      \item  $\teta_i(0)$ reduces modulo $H^{2, 0}$
        to the normalized class $\eta_i(0)$ in $H^{1, 1}(X)_{prim}$.
      \end{enumerate}
      
      Write $\lambda_j(t)=
      {}^t[\int_{c_1}\teta_i(t),\ldots,\int_{c_N}\teta_j]$,
      so that $(\eta_i(t),\eta_j(t))=(\lambda_i(t).\lambda_j(t))$.
      So $(\eta(0),\eta_j(0))$ equals the
      specialization of $\eta_i(t),\eta_j(t)$.
    \end{proof}
  \end{theorem}
  \begin{corollary} Regard $H$ as a multi-valued map
    $H=[\eta_1,\ldots, \eta_N]:\sje^{gen}\to \sT$; then $H$
    satisfies the ecliptic equation.
    \noproof
  \end{corollary}
  
  \begin{remark} Although the Gauss--Manin connection
    does not obviously restrict to a flat holomorphic
    connection on $\sH^{1, 1}_{prim}$ over $\tsje^{gen}$, there is such
    a connection 
    on the pull back of $\sH^{1, 1}_{prim}$ over $\J\E^{gen}$.
    That is, there exists such a connection locally on $\tsje^{gen}$
    and the entries of its connection matrix are the functions $f_{i j}$.
    These locally defined connections are not compatible;
    however, elimination of the unknown coefficients $f_{i j}$ in each
    of the connection matrices gives rise to a non-linear pde
    on each chart and these are
    compatible; they glue to give the ecliptic equation.
  \end{remark}
\end{section}
\begin{section}{Generic immersivity}
  Assume throughout this section that
    $4(1+h-q)\ge 2q+1$.  

  \begin{theorem}\label{imm}
    $H_L:\J\E^{gen}\to \sT_\Lambda$
    is generically immersive.
    \begin{proof} We need three lemmas.
      \begin{lemma} $\sje^{gen}$ dominates $\sM_q$ and
        is irreducible.
        \begin{proof} Suppose $C\in\sM_q$ and
          $L\in\Pic^{1+h-q}_C$. Then $4L$ and $6L$
          are very ample, so that a general member
          of $H^0(C, 4L)\times H^0(C, 6L)$
          defines a surface in $\sje_{h, q}^{gen}$. Since
          the universal $\Pic^{1+h-q}$ over $\sM_q$
          is irreducible the lemma is proved.
        \end{proof}
      \end{lemma}

    \begin{lemma} Suppose that $m\ge 8$.
      Then any transitive subgroup
      of the symmetric group $\mathfrak S_m$ which contains
      a transposition and a copy of
      the alternating group $\mathfrak A_r$
      for some $r > m/2$ is equal to $\mathfrak S_m$.
      \begin{proof} Consideration of the group
        generated by the transpositions in such a subgroup $\G$,
        say,
        shows that there are integers $n, k$ such that
        $n\ge 2$, $m = n k$ and
        $\G$ is a wreath product $\mathfrak S_n\wr K$
        for some transitive subgroup $K$ of $\mathfrak S_k$.
        Since $r > m/2 \ge k$ the resulting
        homomorphism
        ${\mathfrak A}_r\to K$
        is trivial. Then there is an embedding
        ${\mathfrak A}_r\inj \mathfrak S_n$,
        so that $m/k=n\ge r>m/2$ and then $n=m$.
      \end{proof}
    \end{lemma}

    \begin{lemma}    If $\sZ$ is the universal
      ramification divisor then the Galois group of $\sZ\to\sje$ is
      $\mathfrak S_N$.
      \begin{proof} According to Cor. 10.3 of \cite{SB1}
        $\sZ$ is irreducible, because it is
        dominated by the universal
        $\vert 4L\vert\times\vert 6L\vert$,
        so that the Galois group $\G$, say, is a transitive subgroup
        of $\mathfrak S_N$.

        To see that $\G$ contains a transposition, take a surface
        in $\sje_{h, q}$ defined by $\phi:C\to\bsell$ that has
        a point of triple ramification. The corresponding
        $2$-dimensional variation of $X$ that is constructed
        in \cite{SB1} yields a point where $\sZ\to\sje$ is
        simply ramified, and so $\G$ contains a transposition.

        Since 
        $\vert 6L\vert$ is very ample, the 
        Galois group of $\vert 6L\vert$
        is $\mathfrak S_r$, where
        $r=6(h+1-q)-q$.
        Then
        there is an injective homomorphism
        $\mathfrak S_r\to \G$.
        Since $r > N/2$,
        from our assumption, we conclude
        by the previous lemma.
      \end{proof}
    \end{lemma}

    \begin{corollary}\label{refl} $\tsje$ is irreducible.
      \begin{proof} Consider the $\mathfrak S_N$-cover
        $\sje^{gen, *}\to\sje^{gen}$
        of Section \ref{3} and observe that
        $\sje^{gen, *}\to\sje^{gen}$
        is the Galois closure of $\sZ\to\sje$,
        so that $\sje^{gen, *}$ is irreducible. Then to prove
        that $\tsje$ is irreducible it is enough to exhibit
        a monodromy operator that changes the sign of
        exactly one of the $\eta_i$ and leaves the others
        invariant. For this, we essentially repeat
        an argument from the proof of Lemma \ref{monodr}.

        Let the branch point $P_N$
        of $C\to\bsell$ tend to a point where $j=\infty$.
        Then we get a $1$-parameter family
        $\sX\to\Delta$ where the total space $\sX$
        is smooth and the closed fibre $\sX_0$ has an
        $A_1$-singularity. The vanishing cycle $\delta$
        lies in $H^{1, 1}_{prim}$ and is orthogonal
        to $\eta_i$ for $i\ne N$, so that $\delta$
        is a multiple of $\eta_N$.
        The monodromy is of order $2$
        and acts trivially on $\eta_i$ for every $i\ne N$,
        and therefore changes the sign of $\eta_N$.
      \end{proof}
    \end{corollary}

    To prove that $H_L$ is generically immersive
    we can regard $H_L$ as a multi-valued holomorphic map
    from $\tsje^{gen}$ to $\sT$. Observe that $\tsje^{gen}$
    is a cover of $\sje^{gen, *}$ and so is also a cover of $\sZ$.
    
    The Galois group $G_N$ 
    permutes the co-ordinate frame
    $\partial_1, \ldots ,\partial_N$ on $\tsje^{gen}$ and the signed
    vectors $\pm\eta_1,\ldots, \pm\eta_N$. Moreover, the kernel of
    $H_{L *}$ at
    the generic point $\xi$ of $\tsje^{gen}$ is preserved by $G_N$,
    so that if $H_{L}$ is not generically immersive
    then there is a function $f$ on $\bsell$
    such that either
    \begin{enumerate}
    \item $\ker H_{L *}(\xi)$  contains $\sum_1^N f_i \partial_i$ or
    \item $\ker H_{L *}(\xi)$  contains $f_i \partial_i- f_j \partial_j$
      for every $i, j$, \end{enumerate}
    \noindent where $f_i$ is the pull back of $f$ to $\tsje^{gen}$
    under the $i$'th projection $\tsje^{gen}\to \sZ \to \bsell$.

    Suppose first that $\ker H_{L *}(\xi)$
    contains $\sum f_i \partial_i$.
    There is an action of $\GG_m$ on $\bsell$ (it is the integral of
    the vector field $\sum\pd_i$) and so on
    $\sje^{gen}$. If $X$ is given by the Weierstrass equation
    $$Y^2=4X^3-g_4X-g_6$$
    where $g_n \in H^0(C, \phi^*M^{\otimes n})$
    then \cite{SB2} the closure
    of the $\GG_m$-orbit through the point $X$ of $\sje^{gen}$
    is the \emph{Gauss--Eisenstein} pencil
    $$(Y^2=4X^3-\lambda g_4X-\mu g_6)_{(\lambda, \mu) \in \P^1}\ ;$$
    at $\lambda =0$ and at $\mu=0$ the surface has constant
    $j$-invariant and complex multiplication
    by $\Q(\zeta_6)$ and $\Q(\zeta_4)$,
    respectively. Moreover, since this $\GG_m$-action
    is defined over $\Q$, the corresponding vector
    field on $\sje^{gen}$ is invariant under the Galois
    group $G_N$
    and so is proportional to $\sum f_i \partial_i$. Then
    the multi-valued map $H_L$ is constant
    on each Gauss--Eisenstein pencil while the ramification
    divisor $Z$ is constant. Therefore, in each such pencil,
    there is no monodromy
    on the classes $[\eta_i]$. However,
    suppose that $Z$ meets the locus $(j=\infty)$
    in the point $P_i$, so that the surface
    acquires a node. Let $\delta\in H^{1, 1}(X)_{prim}$
    be the corresponding vanishing cycle. Then
    $\eta_i$ and $\delta$ are both orthogonal
    to $\eta_j$ for every $j\ne i$, so that
    $\eta_i$ is proportional to $\delta$. So
    $\eta_i$ is not invariant under
    monodromy; this contradiction completes the
    proof in this case.

    Suppose instead that 
    $\ker H_{L *}(\xi)$ contains
    $f_i \partial_i - f_j \partial_j$ for all $i, j$. Say
    $Z=\{P_1,\ldots ,P_N\}$, $E_i=f^{-1}(P_i)$
    and $\eta_i\in H^0(X,\Omega^2_X(2E_i))_{2k}$.
    Suppose $j\ge 3$. In \cite{SB1} we showed that
    $$\partial_i(\eta_j)=\eta_j(P_i)\eta_i$$
    for $i=1, 2$ when $\eta_j(P_i)$ is the value of the
    $2$-form $\eta_j$ when written in terms of local
    co-ordinates. Since
    $f_1 \partial_1(\eta_j)= f_2 \partial_2(\eta_j)$ it follows
    that $\eta_j$ vanishes at both $P_1$ and $P_2$.
    Therefore $\eta_j$ vanishes at $P_k$ for every $k\ne j$.
    However, $H^0(X,\Omega^2_X(2E_j - \sum_{k\ne j} E_k))=0,$
    and the result is proved.
  \end{proof}
\end{theorem}
\end{section}

\begin{section}{Some formulae}
  The notation will be the same as before.
  \begin{proposition}\label{above} Suppose that $P\in Z$
    and that $z_P$ is a co-ordinate on $C$ at $P$
    such that $z_P^{e_P}$ is a co-ordinate on
    $\bsell$. Assume also that $v_P$ is a
    fibre co-ordinate on $X$ along $E_P$
    that pulls back from a fibre co-ordinate on the universal
    elliptic curve and that vanishes along the zero section
    $C_0$ of $f:X\to C$.
    Suppose that $\eta\in H^0(X,\Omega^2_X(\sum_P e_PE_P))$
    and consider the local expansion
    $$\eta=\sum_{j> -e_P}a_{P, j}z_P^{j-1}dz_P\wedge dv_P$$
    at $P$.
    Then $\eta$ is of the second kind if and only if $a_{P, 0}=0$
    for every $P\in Z$.
    \begin{proof} Drop the subscript $P$,
      so that $v_P=v$, $z_P=z$ and $e_P=e$.
      Then there is a \nbd $\Delta$
      of $P$ in $C$
      such that, if $X_\Delta=f^{-1}(\Delta)$,
      the morphism $X_\Delta\to \Delta$ is defined by
      $$y^2=4x^3-g_4x-g_6$$
      where $g_4,g_6$ are functions of $z^e$
      and $x=\wp(v),\ y=\wp'(v)$.
      Therefore there is an action of
      the group
      $\mu_{e}$ of $e$'th roots of unity
      on $X_\Delta$ defined in terms of the local
      co-ordinates $(z, v)$ by
      $\chi^*(z, v)=(\chi z, v)$ for
      the generator
      $\chi=\exp(2\pi i/e)$ of $\mu_e$.

      Suppose that $\eta$ is of the second kind
      and let $\eta_\Delta$ denote its restriction to $X_\Delta$.
      Then $\g^*\eta_\Delta$ is of the second kind for all
      $\g\in \mu_e$. Take $\g=\chi$; then
      $$\sum_{r=0}^{e-1}\g^{r*}\eta_\Delta=
      e\sum_{k\ge 0}a_{k e}z^{k e}\frac{dz}{z}\wedge dv,$$
      so that $\sum_{r=0}^{e-1}\g^{r*}\eta_\Delta$ is of the second
      kind and has a simple pole. Therefore it is holomorphic,
      so that $a_0=0$.

      Conversely, suppose that $\eta\in H^0(X,\Omega^2_X(eE))$
      and that $\Res_E\eta\ne 0$. Since $\mu_e$ acts trivially on
      $E$, and so on
      $H^1(E,\C)$, it follows that
      $$\Res_E(\g^*\eta_\Delta)=\g^*\Res_E\eta=\Res_E\eta,$$
      so that $\Res_E(\sum_r\g^{r*}\eta_\Delta)=e\Res_E\eta\ne 0$.
      But
      $$\sum_r\g^{r*}\eta_\Delta=
      e\sum_{k\ge 0}a_{k e}z^{k e}\frac{dz}{z}\wedge dv$$
      and so has a simple pole; since its residue
      along $E$ is non-zero, it follows that $a_0\ne 0$.
    \end{proof}
  \end{proposition}
  \begin{proposition} For every $P\in Z$
    we have
    $$\dim H^0(X,\Omega^2_X(e_P P))_{2k}/H^0(X,\Omega^2_X)=e_P -1.$$
    \begin{proof} Recall first that $f^*$ identifies the linear systems
      $\vert K_X\vert$ and $\vert K_C+\phi^*M\vert$.

      The cohomology of the short exact sequence
      $$0\to\sO_C(K_C+\phi^*M)\to\sO_C(K_C+\phi^*M+e_P P)
      \to \sO_{e_P P}\to 0$$
      of sheaves on $C$ shows that
      $$\dim H^0(C,\sO_C(K_C+\phi^*M+e_P P))/
      H^0(C,\sO_C(K_C+\phi^*M)) = e_P.$$
      By the previous proposition there is exactly one linear
      condition on a member of $H^0(X,\Omega^2_X(e_P P))$
      to be of the second kind.
    \end{proof}
  \end{proposition}
  So we can write down a basis of $H^{1, 1}(X)_{prim}$ that is
  represented by $2$-forms of the second kind.
  \begin{theorem} If $P_1, P_2$ are distinct points in $Z$
    and $\eta_{i}\in H^0(X, \Omega^2_X(e_{P_i} P_i))_{2k}$
    for $i=1, 2$ then $(\eta_{1}, \eta_{2})=0$.
    \begin{proof} After pulling back by some finite cover
      $D\to C$ that is unramified over $Z$ we get a surface
      $Y\to D$ defined by $\phi_1:D\to\bsell$ such that
      $\deg\phi_1^*M>1$. Then $\vert K_Y\vert$ has no base
      points so that, by the results of \cite{SB1}, the pull backs
      of the classes $\eta_{i}$ are orthogonal on $Y$.
      Then the $\eta_{i}$ are orthogonal on $X$.
    \end{proof}
  \end{theorem}
\end{section}
\begin{section}{Infinitesimal Torelli with base points}
  As already mentioned, Saito \cite{S} (see also
  \cite{EGW}, \cite{I} and \cite{Kl})
  proved infinitesimal Torelli when $\Bs\vert K_X\vert$, which
  we identify with $\Bs\vert K_C+\phi^*M\vert$,
  is empty. Since $\deg \phi^*M>0$, if $\Bs\vert K_X\vert$ is not
  empty then it consists of a single point $Q\in C$,
  where $\phi^*M=\sO_C(Q)$,
  and also $h=q$.

  \begin{theorem} Assume that 
    $\Bs\vert K_X\vert =\{Q\}$. Then
    infinitesimal Torelli holds for $X$ if and only
    if $Q\not\in Z$. If $Q\in Z$ then the derivative
    of the period map is of corank $1$.
    \begin{proof} Suppose that $P\in Z$ and put $e_P=e$.
      Then consider the $1$-parameter variation
      $\sX\to\Delta_t$ of $X$
      constructed from the glueing
      $$z^e = w^e + e t w^{e-2};$$
      this is a sub-variation of the $(e-1)$-parameter
      variation detailed in \cite{SB1}. Observe that
      $$w= z(1-t z^{-2})$$
      modulo $t^2$. Consider
      $\omega^{(j)}(t)\in H^0(\sX_t, \Omega^2_{\sX_t})$ normalized
      by the requirement that
      $\int_{A_i}\omega^{(j)}(t)=\delta_i^j$
      for appropriate $2$-cycles $A_i$ on $X$ (or on $\sX_t$) supported
      away from the fibre $E_Q$. Expand $\omega^{(j)}(t)$ as
      $$\omega^{(j)}(t)=\sum_{p, q\ge 0} b^{(j)}_{p q} w^p t^q dw \wedge dv$$
      for some fibre co-ordinate $v$. Substituting
      $w= z(1-t z^{-2})$ gives
      $$\omega^{(j)}(t)=
      \left(\sum b^{(j)}_{p q}z^pt^q-
        \sum b^{(j)}_{p q}(p-1)z^{p-2}t^{q+1}\right)dz\wedge dv$$
      modulo $t^2$. Write $\omega^{(j)}(t)=\omega^{(j)}+t\eta^{(j)};$
      then $\omega^{(j)}=\sum b_{p 0}^{(j)}z^p dz\wedge dv$
      and $\eta^{(j)}= b^{(j)}_{0 0}z^{-2}dz\wedge dv$ modulo
      $H^0(X,\Omega^2_X)$.

      Now $b^{(j)}_{00}$ is the constant
      coefficient of $\omega^{(j)}$ and so vanishes
      if $P=Q$. Suppose indeed that $P=Q$; then
      the class of $\eta^{(j)}$ in $H^{1, 1}_{prim}(X)$ also vanishes
      for every $j$, so that ${\partial/\partial t}:
      H^{2, 0}(X)\to H^{1, 1}_{prim}(X)$ is the zero map
      and indeed infinitesimal Torelli fails.
      
      On the other hand, if $Q\not\in Z$ then
      for each $P\in Z$ the image of $\partial_{t_P}$
      is the line $\C\eta_P$ and it follows
      that no linear combination of the
      $\partial_{t_P}$
      can vanish, so that infinitesimal Torelli holds in this case.
    \end{proof}
  \end{theorem}
\end{section}
\begin{section}{Intersecting forms with curves}
  We give a formula
  for the cup product $[\eta]\cup [D]$ when
  $\eta$ is a $2$-form of the second kind
  and $D$ is a curve in $X$ that is unramified
  over $C$ at the points of the ramification divisor $Z$.

  Recall that $X$ is provided
  with a zero section $C_0$ as part of its data.
  
  \begin{theorem}
    Suppose that $\eta\in H^0(X,\Omega^2_X(\sum_{P\in Z}
    e_PE_P))_{2k}$
    and $D$ is as above. Assume that
    we have chosen local co-ordinates
    as in Proposition \ref{above}. Assume further that
    $D$ is defined locally at each point $Q$
    of $D\cap E_P$ by
    $$v_P=v_P(z_P)=\sum_{m\ge 0} c_{Q, m}z_P^m$$
    and that the corresponding local expansion of $\eta$ is
    $$\eta=\sum_{0\ne n> -e_P}a_{P, n}z_P^{n-1}dz_P\wedge dv_P$$
    for every $P\in Z$.
    Then
    $$[\eta]\cup_{dR} [D]=
    -\sum_{P\in Z}\sum_{m= 1}^{e_P-1}a_{P, -m}
    \sum_{Q\in D\cap E_P}c_{P, m}.$$
    \begin{proof} Without loss of generality 
      we can assume that the polar divisor
      $(\eta)_\infty$ is supported on a single fibre $E_P$.

      Choose small open
      discs $P\in\Delta'\subsetneq \Delta\subsetneq C$
      and a small tubular \nbd $X^0$ of $D$ in $X$.
      Set $R=\bDelta-\Delta'$, a closed annulus,
      and let $\g=\partial\bDelta,\g'=\partial\bDelta'$.
      Then there is a $C^\infty$ $(1,0)$-form $\Psi$
      on $X^0$ such that
      \begin{enumerate}
      \item $\Psi=0$ on $f^{-1}(C-\Delta)\cap X^0$ and
      \item $\Psi= \sum_{n\ne 0}\frac{a_n}{n}z^n dv$
        on $f^{-1}(\Delta')\cap X^0$, so that
        $d\Psi=\eta$ there.
      \end{enumerate}
      Put $\txi=\eta-d\Psi$, a $C^\infty$ $2$-form on $X^0$.
      Then $[\eta]\vert_{X^0}=[\txi]$, $\txi$ is zero on
      $f^{-1}(\Delta')\cap X^0$ and $\txi$ is holomorphic
      on $f^{-1}(C-\Delta)\cap X^0$. So
      $$[\eta]\cup_{dR} [D] =  [\eta]\vert_{X^0}\cup_{dR} [D]
      =\frac{1}{2\pi i}\int_D\txi\vert_D =
      \frac{1}{2\pi i}\int_{D\cap f^{-1}(R)}\txi\vert_D,$$
      since $\txi\vert_D$ is zero on $D\cap f^{-1}(\Delta')$
      and is a holomorphic $2$-form on $D\cap f^{-1}(C-\Delta)$,
      so is zero there. So
      \begin{eqnarray*}
        2\pi i\ [\eta]\cup_{dR} [D] &=& \int_{D\cap f^{-1}(R)}(\eta-d\Psi)\vert_D\\
        &=&\int_{D\cap f^{-1}(R)}\eta\vert_D-\int_{D\cap
                           f^{-1}(R)}d\Psi\vert_D\\
                       &=& 0-\int_{D\cap f^{-1}(R)}d\Psi\vert_D\\
                       &=& \int_{D\cap f^{-1}(\g')}\Psi\vert_D-
                           \int_{D\cap f^{-1}(\g)}\Psi\vert_D\
                           \textrm{(Stokes)}\\
                       &=& \int_{D\cap f^{-1}(\g')}\Psi\vert_D\
                           \textrm{since}\ \Psi =0\
                           \textrm{on}\ f^{-1}(C-\Delta)\\
                                   &=& 2\pi i\ \sum_{Q\in D\cap E_P}
                                       \Res_Q\Psi\vert_D\ \textrm{(Cauchy)}.
      \end{eqnarray*}
      Now
      $$\Psi=\sum_{n\ne 0}\frac{a_n}{n}z^n dv=
      \sum_{n\ne 0,\ m\ge 0}\frac{m}{n}c_{Q, m} a_n z^{m+n-1}dz,$$
      so that
      $$ [\eta]\cup_{dR} [D]=\sum\Res_Q\Psi
      =\sum_Q\sum_{m\ge 1}\frac{m}{-m}c_{Q, m}a_{-m}
      =-\sum_{m=1}^{e_P-1}a_{-m}\sum_Qc_{Q, m}.$$
    \end{proof}
  \end{theorem}

  \begin{corollary}\label{4.2} If $f:X\to C$ is written in
    Weierstrass form as
    $$y^2=4x^3-g_4x-g_6,$$
    if $\phi:C\to \bsell$ is simply ramified
    and if every $v_P$ is chosen so that
    $dv_P=dx/y$ then
    $$[\eta]\cup_{dR}[D]=-\sum_{P\in Z}
    a_{P, -1}\frac{dx}{dz_P}(P)\sum_{Q\in D\cap E_P}
    \frac{1}{y(Q)}.$$
    \begin{proof} Put $z_P=z$ and $v_P=v$. Then,
      in terms of a period $\tau=\tau(z)$ of $E_z$
      we have 
      $$x=\wp(v,\tau),\ y=\frac{\partial\wp(v,\tau)} {\partial v},$$
      so that
      $$\frac{d x}{d z} = y\frac{d v}{d z}
      +\frac{\partial\wp}{d\tau}\frac{d\tau}{d z}.$$
      Now $\frac{d\tau}{d z}\vert_{z=0}=0$, since
      $\phi$ is ramified at $P$, so that
      $$c_{Q,1}=
      \frac{d v}{d z}(Q)
      =\frac{1}{y(Q)}\frac{d x}{d z}(P).$$
    \end{proof}
  \end{corollary}

  \begin{remark}
    \part[i] If the ramification is not simple then
    we can calculate the coefficients
    $c_{P, n}$ in a similar way via Fa{\`a} di Bruno's formula.
    We omit the details, however.


    \part[ii] In particular, $[\eta]\cup [C_0]=0,$
    where $C_0$ is the zero section. Since also
    $[\eta]\cup [f]=0$ for a fibre $f$, it
    follows that $[\eta]$ lies in $H^2_{prim}(X,\C)$,
    and so in $\Fil^1H^2_{prim}(X,\C)$.
  \end{remark}
\end{section}
\begin{section}{Rational surfaces}
  For rational surfaces the multi-valued map $H_L$
  is a single-valued map $H:\tsje^{gen}\to\sT$.
  We say something about how to calculate it.
  
  Suppose that $f:X\to C =\P^1$ is a generic
  Jacobian elliptic surface which is rational,
  that $C_0$ is the given zero section
  and that $\phi$ is a fibre.
  Then there are eight further sections $C_1,\ldots ,C_8,$
  disjoint from each other and from $C_0$, and
  the curves
  $C_0,C_1,\ldots ,C_8$ are the exceptional
  curves of a birational contraction $\pi:X\to\P^2$.
  The lattice $H^2_{prim}(X,\Z)$ is isomorphic
  to the root lattice $E_8(-1)$ and if
  $h$ is the class of a line in $\P^2$ and
  $\a_1,\ldots ,\a_8$ are defined by
  $$\a_1=C_1-C_2,\ \a_2=2h-C_1-C_2-C_3\ \textrm{while}\
  \a_j=C_{j-1}-C_j\ \textrm{for}\ j\ge 3$$
  then $(\a_1,\ldots ,\a_8)$ is a root basis of
  $H^2_{prim}(X,\Z)$. Note that
  $\phi=3h-\sum_0^8 C_i.$

  We can calculate each intersection number
  $[\eta_i]\cup_{dR} [C_j]$, and therefore each number
  $[\eta_i]\cup_{dR} [\a_j]$. Then the period matrix $H(X)$
  is
  the $8\times 8$ matrix
  $$H(X)=\left(
    [\eta_i]\cup_{dR} [\a_j]
  \right)
  $$
  and we can calculate it as follows.
  \begin{enumerate}
    
  \item Take a field
  $K$ of characteristic zero and take eight points
  $Q_1,\ldots ,Q_8\in\P^2(K)$. Let $h$ denote the class
  of a line in $\P^2$.

\item Check that these points are in general position
  in the usual sense: they are distinct; no three are collinear;
  no six lie on a conic; there is no cubic through all of them
  that is singular at one of them.
  
\item Calculate a basis $(U_1,U_2)$ of
  $H^0(\P^2_K,\sO(3h-\sum Q_i))$.
\item Construct a basis $(U_1^2,U_1U_2,U_2^2,V)$
  of $H^0(\P^2_K,\sO(6h-2\sum Q_i))$.

\item Construct a basis
    $(U_1^3,U_1^2U_2,\ldots ,VU_2,W)$
    of $H^0(\P^2_K,\sO(9h-3\sum Q_i))$.
  \item The forms $(U_1,U_2,V,W)$ define a
    rational map
    $$\rho:\P^2-\to \P(1,1,2,3)=\Proj K[U_1,U_2,V,W]$$
    whose image is a surface $\Sigma$.
    This is the anticanonical model of
    the smooth del Pezzo surface $\Bl_{Q)1,\ldots ,Q_8}\P^2$
    of degree $1$ and is embedded as a sextic surface in $\P(1,1,2,3)$.
  \item After calculating the images $\rho(x)$ for
    sufficiently many $K$-points $x\in\P^2$ (or we could
    calculate $\rho(\eta)$ where $\eta$ is the generic point)
    we can interpolate through those points to
    find the defining sextic equation $F=0$ of $\Sigma$.
  \item Make an explicit change of co-ordinates
    so that
    $F=W^2-(4V^3-g_4V-g_6)$
    where $g_n=g_n(U_1,U_2)$ is homogeneous of degree $n$.
  \item Let $Q_9$ denote the ninth base point
    of the pencil $\vert 3h-\sum_1^8 Q_i\vert$
    and put $X=\Bl_{Q_9}\Sigma$. Let
    $C_0$ denote the exceptional curve of the blow-up;
    then $C_0$ is a section of the elliptic
    fibration $f:X\to\P^1$ defined by
    $f(Q)=(U_1(Q),U_2(Q))$.
  \item If $j\ne k$ then the line $D_{jk}=h-Q_j-Q_k$ is a section
    of $f$. The isomorphism $D_{i j}\to\P^1$
    induced by $f$ can be calculated.
  \item The ramification locus $Z$ in $\P^1$
    is defined by the vanishing of the Jacobian covariant
    $$J=\det\left[
      \begin{array}{cc}
        {{\partial g_4}/{\partial U_1}}&{{\partial g_4}/{\partial U_2}}\\
        {{\partial g_6}/{\partial U_1}}&{{\partial g_6}/{\partial U_2}}
      \end{array}
    \right].
    $$
  \item Check the genericity of the surface $f:X\to\P^1$
    by calculating the discriminant $D=g_4^3-27g_6^2$
    and then checking that the discriminants of $D$ and $J$
    and the resultant of $D$ and $J$ are all non-zero.
    (Of course, this could be carried out by reduction
    modulo $p$.)
  \item Now $K(Z)$ is a separable octic extension of $K$
    and $e_P=2$ for all geometric points $P_i\in Z$.
    Calculate each form $\eta_i=\eta_{P_i}$,
    up to a scalar $\mu_i$,
    by writing down sections of the two-dimensional vector
    space $H^0(X,\Omega^2_X(2P_i))$ and using Proposition \ref{above}.
  \item Calculate each intersection number $[\eta_i]\cup_{dR} [D_{jk}]$
    via Corollary \ref{4.2} and then calculate each
    $[\eta_i]\cup \a_l$ from the description of $\a_l$ as a
    $\Q$-linear combination of the $D_{jk}$.
  \item Calculate each $[\eta_i]$ as a linear combination
    of the $\a_l$ from the $8\times 8$ matrix
    $([\eta_i]\cup \a_l)$.
  \item We know that the $\eta_i$ are orthogonal; we can
    now make them orthonormal. That is, we can determine the
    scalars $\mu_i$.
  \item We have now calculated $H(X)$ up to
    a single sign if $K(Z)$ is a field.
  \end{enumerate}
\end{section}

\end{document}